\documentclass[11pt]{amsart}
\usepackage{latexsym,amsfonts,amssymb}
\theoremstyle{definition}

\theoremstyle{remark}

\numberwithin{equation}{section}

\begin{document}
\baselineskip 0.55cm

\noindent  \vspace{0.5in}

\title{An extremal problem
on potentially  $K_{r+1}-(kP_2\cup tK_2)$-graphic sequences}
\author{Chunhui Lai}
\author{Yuzhen Sun}

\address{Chunhui Lai(correspondent author): Department of Mathematics,
Zhangzhou Teachers College,  Zhangzhou 363000, P. R. China}
\email{zjlaichu@public.zzptt.fj.cn, laich@winmail.cn}
\address{Yuzhen Sun: Department of Information technology, Zhangzhou Institute of
Education, Zhangzhou, Fujian 363000, P. R. of CHINA.}
\thanks{Rearch is supported by NNSF of China(10271105)
and by NSF of Fujian(Z0511034), Fujian Provincial Training
Foundation for "Bai-Quan-Wan Talents Engineering", Project of Fujian
Education Department, Project of Zhangzhou Teachers College  }



\begin{abstract} A graphic sequence $S$ is potentially $K_{m}-H$-graphical if it has
a realization containing a $K_{m}-H$ as a subgraph. Let
$\sigma(K_{m}-H, n)$ denote the smallest degree sum such that every
$n$-term graphical sequence $S$ with $\sigma(S)\geq \sigma(K_{m}-H,
n)$ is potentially $K_{m}-H$-graphical.  In this paper, we determine
$\sigma (K_{r+1}-(kP_2\bigcup tK_2), n)$ for
    $n\geq 4r+10, r+1 \geq 3k+2t,
    k+t \geq 2,k \geq 1, t \geq 0$ . To now, the problem
of determined $\sigma (K_{r+1}-H, n)$ for $H$ not containing a cycle
on 3 vertices and sufficiently large $n$ has been solved.
\end{abstract}
\maketitle

\par
 \section{Introduction}
\par

  The set of all non-increasing nonnegative integer sequences $\pi=$
  ($d_1 ,$ $d_2 ,$ $...,$ $d_n $) is denoted by $NS_n$.
  A sequence
$\pi\epsilon NS_n$ is said to be graphic if it is the degree
sequence of a simple graph $G$ on $n$ vertices, and such a graph $G$
is called a realization of $\pi$. The set of all graphic sequences
in $NS_n$ is denoted by $GS_n$. A graphical sequence $\pi$ is
potentially $H$-graphical if there is a realization of $\pi$
containing $H$ as a subgraph, while $\pi$ is forcibly $H$-graphical
if every realization of $\pi$ contains $H$ as a subgraph. If $\pi$
has a realization in which the $r+1$ vertices of  largest degree
induce a clique, then $\pi$ is said to be potentially
$A_{r+1}$-graphic. Let $\sigma(\pi)=d(v_1 )+d(v_2 )+... +d(v_n ),$
and $[x]$ denote the largest integer less than or equal to $x$. If
$G$ and $G_1$ are graphs, then $G\cup G_1$ is the disjoint union of
$G$ and $G_1$. If $G = G_1$, we abbreviate $G\cup G_1$ as $2G$. We
denote $G+H$ as the graph with $V(G+H)=V(G)\bigcup V(H)$ and
$E(G+H)=E(G)\bigcup E(H)\bigcup \{xy: x\in V(G) , y \in V(H) \}. $
Let $K_k$, $C_k$, $T_k$, and $P_{k}$ denote a complete graph on $k$
vertices,  a cycle on $k$ vertices, a tree on $k+1$ vertices, and a
path on $k+1$ vertices, respectively. Let $K_{m}-H$ be the graph
obtained from $K_{m}$ by removing the edge set $E(H)$ of the graph
$H$ ($H$ is a subgraph of $K_{m}$).
\par

Given a graph $H$, what is the maximum number of edges of a graph
with $n$ vertices not containing $H$ as a subgraph? This number is
denoted $ex(n,H)$, and is known as the Tur\'{a}n number. This
problem was proposed for $H = C_4$ by Erd\"os [2] in 1938 and in
general by Tur\'{a}n [20]. In terms of graphic sequences, the number
$2ex(n,H)+2$ is the minimum even integer $l$ such that every
$n$-term graphical sequence $\pi$ with $\sigma (\pi) \geq l $ is
forcibly $H$-graphical. Here we consider the following variation:
determine the minimum even integer $l$ such that every $n$-term
graphical sequence $\pi$ with $\sigma(\pi)\ge l$ is potentially
$H$-graphical. We denote this minimum $l$ by $\sigma(H, n)$.
Erd\"os,\ Jacobson and Lehel [4] showed that $\sigma(K_k, n)\ge
(k-2)(2n-k+1)+2$ and conjectured that equality holds. They proved
that if $\pi$ does not contain zero terms, this conjecture is true
for $k=3,\ n\ge 6$. The conjecture is confirmed in
[5],[15],[16],[17] and [18]. Recently, Ferrara, Gould and Schmitt
[7] proved the conjecture  using graph theoretic techniques.
Ferrara, Gould and Schmitt [8] determined in  $\sigma(F_k,n)$ where
$F_k$ denotes the graph of $k$ triangles intersecting at exactly one
common vertex.
 \par
 Gould,\ Jacobson and
Lehel [5] also proved that  $\sigma(pK_2, n)=(p-1)(2n-2)+2$ for
$p\ge 2$; $\sigma(C_4, n)=2[{{3n-1}\over 2}]$ for $n\ge 4$. Luo [19]
characterized the potentially $C_{k}$ graphic sequences for
$k=3,4,5.$  Gupta, Joshi and Tripathi [6] gave a necessary and
sufficient condition for the existence of a tree of order
 $n$ with a given degree set. Meng-Xiao Yin, Jian-Hua Yin [28] characterize the
potentially $(K_5-e)$-positive graphic sequences and give two simple
necessary and sufficient conditions for a positive graphic sequence
$\pi$ to be potentially $K_5$-graphic, where $K_r$ is a complete
graph on $r$ vertices and $K_r-e$ is a graph obtained from $K_r$ by
deleting one edge. Moreover, they also give a simple necessary and
sufficient condition for a positive graphic sequence $\pi$ to be
potentially $K_6$-graphic. Gould et al. [5] determined
$\sigma(K_{2,2},n)$ for $n \geq 4$. Yin et al. [23-26] determined
$\sigma(K_{r,s},n)$ for $s\geq r\geq 2$ and sufficiently large $n$.
Lai [10] determined  $\sigma (K_4-e, n)$ for $n\ge 4$.\ Yin,Li and
Mao[22] determined $\sigma(K_{r+1}-e,n)$ for $r\geq 3,$ $r+1\leq n
\leq 2r$ and $\sigma(K_5-e,n)$ for $n\geq5$. Yin and Li[21] gave a
good method (Yin-Li method) of determining the values
$\sigma(K_{r+1}-e,n)$ for $r\geq2$ and $n\geq3r^2-r-1$ (In fact, Yin
and Li[21] also determining the values $\sigma(K_{r+1}-ke,n)$ for
$r\geq2$ and $n\geq3r^2-r-1$). After reading[21],  using Yin-Li
method Yin [27] determined $\sigma (K_{r+1}-K_{3}, n)$ for
    $n\geq 3r+5, r\geq 3$. Lai $[11]$ determined
$\sigma(K_{5}-K_{3}, n)$ for $n\geq 5$ .
      Lai [12,13] determined
    $\sigma (K_{5}-C_{4}, n),\sigma (K_{5}-P_{3}, n)$ and
    $\sigma (K_{5}-P_{4}, n),$ for  $n\geq 5$. Determining $\sigma(K_{r+1}-H,n)$, where $H$
    is a tree on 4 vertices is more useful than a cycle on 4
    vertices (for example, $C_4 \not\subset C_i$, but $P_3 \subset C_i$ for $i\geq 5$).
    So, after reading[21] and [27],  using Yin-Li method Lai and Hu[14]
    determined  $\sigma (K_{r+1}-H, n)$ for
    $n\geq 4r+10, r\geq 3, r+1 \geq k \geq 4$ and $H$ be a graph on $k$
    vertices which
    containing a tree on $4$ vertices but
     not containing a cycle on $3$ vertices and $\sigma (K_{r+1}-P_2, n)$ for
    $n\geq 4r+8, r\geq 3$.
    In this paper,  using Yin-Li method we prove the following two
theorems.\par

{\bf  Theorem 1.1.} If $r\geq 4$ and $n\geq 4r+10$, then $\sigma
(K_{r+1}-(P_2\bigcup K_2), n)= (r-1)(2n-r)-2(n-r).$

\par
{\bf  Theorem 1.2.} If $n\geq 4r+10,  r+1 \geq 3k+2t,
    k+t \geq 2,k \geq 1, t \geq 0$, then $\sigma (K_{r+1}-(kP_2\bigcup tK_2), n)=
(r-1)(2n-r)-2(n-r)$.
       \par
       To now, the problem
of determined $\sigma (K_{r+1}-H, n)$ for $H$ not containing a cycle
on 3 vertices and sufficiently large $n$ has been solved.
\par

\section{Preparations}\par
  In order to prove our main result,we need the following notations
  and results.\par
  Let $\pi=(d_1,\cdots,d_n)\epsilon NS_n,1\leq k\leq n$. Let \par
    $$ \pi_k^{\prime\prime}=\left\{
    \begin{array}{ll}(d_1-1,\cdots,d_{k-1}-1,d_{k+1}-1,
    \cdots,d_{d_k+1}-1,d_{d_k+2},\cdots,d_n), \\ \mbox{ if $d_k\geq k,$}\\
    (d_1-1,\cdots,d_{d_k}-1,d_{d_k+1},\cdots,d_{k-1},d_{k+1},\cdots,d_n),
     \\ \mbox{if $d_k < k.$} \end{array} \right. $$
  Denote
  $\pi_k^\prime=(d_1^\prime,d_2^\prime,\cdots,d_{n-1}^\prime)$,where
  $d_1^\prime\geq d_2^\prime\geq\cdots\geq d_{n-1}^\prime$ is a
  rearrangement of the $n-1$ terms of $\pi_k^{\prime\prime}$. Then
  $\pi_k^{\prime}$ is called the residual sequence obtained by
  laying off $d_k$ from $\pi$.\par
    {\bf Theorem 2.1[21]} Let $n\geq r+1$ and $\pi=(d_1,d_2,\cdots,d_n)\epsilon
    GS_n$ with $d_{r+1}\geq r$. If $d_i\geq 2r-i$ for
    $i=1,2,\cdots,r-1$, then $\pi$ is potentially $A_{r+1}$-graphic.
    \par
    {\bf Theorem 2.2[21]} Let $n\geq 2r+2$ and $\pi=(d_1,d_2,\cdots,d_n)\epsilon
    GS_n$ with $d_{r+1}\geq r$. If $d_{2r+2}\geq r-1$ , then $\pi$ is
    potentially $A_{r+1}$-graphic.
    \par
    {\bf Theorem 2.3[21]} Let $n\geq r+1$ and $\pi=(d_1,d_2,\cdots,d_n)\epsilon
    GS_n$ with $d_{r+1}\geq r-1$. If $d_i\geq 2r-i$ for
    $i=1,2,\cdots,r-1$, then $\pi$ is potentially $K_{r+1}-e$-graphic.
    \par
    {\bf Theorem 2.4[21]} Let $n\geq 2r+2$ and $\pi=(d_1,d_2,\cdots,d_n)\epsilon
    GS_n$ with $d_{r-1}\geq r$. If $d_{2r+2}\geq r-1$ , then $\pi$ is
    potentially $K_{r+1}-e$
    -graphic.
    \par
    {\bf Theorem 2.5[9]} Let $\pi=(d_1,\cdots,d_n)\epsilon NS_n$ and $1\leq k\leq
    n$. Then $\pi\epsilon GS_n$ if and only if  $\pi_k^\prime\epsilon
    GS_{n-1}$.
    \par
     {\bf Theorem 2.6[3]} Let $\pi=(d_1,\cdots,d_n)\epsilon NS_n$
     with even $\sigma(\pi)$. Then $\pi\epsilon GS_n$ if and only if
     for any $t$,$1\leq t\leq n-1$,
     $$\sum_{i=1}^t d_i\leq t(t-1)+\sum_{j=t+1}^n
     min \{t,d_j \}.$$
     \par
     {\bf Theorem 2.7[5]} If $\pi=(d_1,d_2,\cdots,d_n)$ is a graphic
     sequence with a realization $G$ containing $H$ as a subgraph, then
     there exists a realization $G^\prime$ of $\pi$ containing H as a
     subgraph so that the vertices of $H$ have the largest degrees of
     $\pi$.
     \par
     {\bf Lemma 2.1 [27]} If $\pi=(d_1,d_2,\cdots,d_n)\epsilon NS_n$ is potentially
     $K_{r+1}-e$-graphic, then there is a realization $G$ of $\pi$
     containing $K_{r+1}-e$ with the $r+1$ vertices
     $v_1,\cdots,v_{r+1}$ such that $d_G(v_i)=d_i$ for
     $i=1,2,\cdots,r+1$ and $e=v_rv_{r+1}$.
     \par

     {\bf Lemma 2.2  [14]} Let $n\geq r+1$ and $\pi=(d_1,d_2,\cdots,d_n)\epsilon GS_n$
with $d_r\geq r-1$ and $d_{r+1}\geq r-2$. If $d_i\geq 2r-i$ for
$i=1,2,\cdots,r-2,$ then $\pi$ is potentially $K_{r+1}-P_2$-graphic.

    \par

{\bf Lemma 2.3  [14]} Let $n\geq 2r+2$ and
$\pi=(d_1,d_2,\cdots,d_n)\epsilon GS_n$ with  $d_{r-2}\geq r$. If
$d_{2r+2}\geq r-1$, then $\pi$ is potentially $K_{r+1}-P_2$-graphic.
\par

\section{ Proof of Main results.} \par

    {\bf  Lemma 3.1.} If $n\geq r+1, r+1 \geq 3k+2t,
    k+t \geq 2,k \geq 1, t \geq 0$,
     then $\sigma (K_{r+1}-(kP_2\bigcup tK_2), n)\geq (r-1)(2n-r)-2(n-r)$.
\par
{\bf Proof.}   Let
$$G =K_{r-2}+  \overline{K_{n-r+2}}$$
Then $G$ is a unique realization of $((n-1)^{r-2}, (r-2)^{n-r+2})$
and $G$ clearly does not contain $K_{r+1}-(kP_2\bigcup tK_2),$ where
the symbol $x^y$ means $x$ repeats $y$ times in the sequence. Thus
$\sigma (K_{r+1}-(kP_2\bigcup tK_2), n) \geq
(r-2)(n-1)+(r-2)(n-r+2)+2= (r-1)(2n-r)-2(n-r).$
\par

 {\bf The Proof of Theorem 1.1 }
       According to Lemma 3.1, it is enough to verify that for
        $n\geq 4r+10$,
       $$\sigma(K_{r+1}-(P_2\bigcup K_2),n)\leq(r-1)(2n-r)-2(n-r).$$
    We now prove that if $n\geq 4r+10$ and $\pi=(d_1,d_2,\cdots,d_n)\epsilon GS_n$
with
$$\sigma(\pi)\geq(r-1)(2n-r)-2(n-r),$$
then $\pi$ is potentially
  $K_{r+1}-(P_2\bigcup K_2)$-graphic.
  \par
   If $d_{r-2}\leq r-1$.
   \par
   (1)Suppose $d_{r-2}= r-1$ and
 $\sigma(\pi)=(r-3)(n-1)+(r-1)(n-r+3)$, then
 $\pi=((n-1)^{r-3},(r-1)^{n-r+3})$. Obviously $\pi$ is potentially
 $K_{r+1}-(P_2\bigcup K_2)$ graphic.
 \par
 (2)Suppose $d_{r-2}= r-1$ and
 $\sigma(\pi)<(r-3)(n-1)+(r-1)(n-r+3)$,
 then
 $$ \begin{array}{rcl}
   \sigma(\pi)&<&(r-3)(n-1)+(r-1)(n-r+3)\\
   &=&   (r-1)(n-1)-2(n-1)+(r-1)(n-r+3)\\
   &=& (r-1)(2n-r)-2(n-r),

   \end{array} $$
 which is
 a contradiction.
 \par
 (3)Suppose $d_{r-2}< r-1$, then
   $$ \begin{array}{rcl}
   \sigma(\pi)&< &(r-3)(n-1)+(r-1)(n-r+3)\\
   &=&   (r-1)(n-1)-2(n-1)+(r-1)(n-r+3)\\
   &=& (r-1)(2n-r)-2(n-r),
   \end{array}
    $$
   which is a
    contradiction.
    \par
    Thus, $d_{r-2}\geq r$ or $\pi$ is potentially
 $K_{r+1}-(P_2\bigcup K_2)$ graphic.
    \par
    If $d_r\leq r-2$, then
    $$ \begin{array}{rcl}
    \sigma(\pi)&=&\sum_{i=1}^{r-1}d_i+\sum_{i=r}^n d_i\\
    & \leq & (r-1)(r-2)+\sum_{i=r}^{n}min \{r-1,d_i \}+\sum_{i=r}^n
    d_i\\
    &=&(r-1)(r-2)+2\sum_{i=r}^n
    d_i\\
    &\leq &(r-1)(r-2)+2(n-r+1)(r-2)\\
    &=& (r-1)(2n-r)-2(n-r)-2\\
    &<& (r-1)(2n-r)-2(n-r),
    \end{array} $$
    which is a
    contradiction. Hence $d_r\geq r-1$.
\par
    If $d_{r+1}\leq r-3.$
    \par
  (1)Suppose $d_r=n-1$, then $d_1 \geq d_2\geq \cdots \geq d_{r-1}\geq
 d_r=n-1$, therefore $d_1=d_2=\cdots=d_r=n-1$. Therefore $d_{r+1} \geq r,$
 which is a contradiction.
 \par
 (2)Suppose $d_r\leq n-2,$ then
    $$\begin{array}{rcl}
    \sigma(\pi)&=&\sum_{i=1}^{r-1}d_i+d_r+\sum_{i=r+1}^n d_i \\
    & \leq &(r-1)(r-2)+
    \sum_{i=r}^{n}min \{r-1,d_i \}+d_r+\sum_{i=r+1}^n d_i \\
    &=&(r-1)(r-2)+min \{r-1,d_r \}+d_r+2\sum_{i=r+1}^n
    d_i \\
    & \leq & (r-1)(r-2)+2d_r+2\sum_{i=r+1}^nd_i \\

    & \leq &(r-1)(r-2)+2(n-2)+2(n-r)(r-3) \\
    &=&(r-1)(2n-r)-2(n-r)-2 \\
    &<&(r-1)(2n-r)-2(n-r), \end{array} $$
     which is a
    contradiction.
    \par
    Thus $d_{r+1}\geq r-2$.
    \par
      If $d_i\geq 2r-i$ for $i=1,2,\cdots,r-2$ or $d_{2r+2}\geq r-1$, then
      $\pi$ is potentially
 $K_{r+1}-(P_2\bigcup K_2)$ graphic($\pi=((n-1)^{r-3},(r-1)^{n-r+3})$) or $\pi$
      is potentially
  $K_{r+1}-P_2$-graphic  by Lemma 2.2 or Lemma 2.3 . Therefore, $\pi$
      is potentially
  $K_{r+1}-(P_2\bigcup K_2)$-graphic. If $d_{2r+2}\leq
  r-2$ and there exists an integer $i$, $1\leq i\leq r-2$ such that
  $d_i\leq 2r-i-1$, then
  $$\begin{array}{rcl}
  \sigma(\pi) &\leq &(i-1)(n-1)+(2r+1-i+1)(2r-i-1)\\
  &&+
  (r-2)(n+1-2r-2)\\
  &=& i^2+i(n-4r-2)-(n-1)\\
  &&+(2r-1)(2r+2)+(r-2)(n-2r-1).
  \end{array} $$
  Since $n\geq 4r+10$, it is  easy to see that
  $i^2+i(n-4r-2)$, consider as a function of $i$, attain its maximum
  value when $i=r-2$. Therefore,
  $$\begin{array}{rcl}
  \sigma(\pi) &\leq &
  (r-2)^2+(n-4r-2)(r-2)-(n-1)\\
  &&+(2r-1)(2r+2)+(r-2)(n-2r-1)\\
  &=&(r-1)(2n-r)-2(n-r)-n+4r+9 \\
  &<&\sigma(\pi),
  \end{array} $$
  which is a
    contradiction.
    \par
    Thus, $\sigma
(K_{r+1}-(P_2\bigcup K_2), n) \leq(r-1)(2n-r)-2(n-r)$ for
    $n\geq 4r+10$.
\par

 {\bf The Proof of Theorem 1.2 } By Lemma 3.1,
 for $n\geq 4r+10,  r+1 \geq 3k+2t,
    k+t \geq 2,k \geq 1, t \geq 0$,
$\sigma (K_{r+1}-(kP_2\bigcup tK_2), n)\geq (r-1)(2n-r)-2(n-r).$
Obviously, for $n\geq 4r+10, r+1 \geq 3k+2t,
    k+t \geq 2,k \geq 1, t \geq 0$,
$\sigma (K_{r+1}-(kP_2\bigcup tK_2), n)\leq
\sigma(K_{r+1}-(P_2\bigcup K_2), n).$ By theorem 1.1, for $n\geq
4r+10, r\geq 4$, $\sigma (K_{r+1}-(P_2\bigcup K_2), n)=
(r-1)(2n-r)-2(n-r).$  Then $\sigma (K_{r+1}-(kP_2\bigcup tK_2), n)=
(r-1)(2n-r)-2(n-r),$  for  $n\geq 4r+10,  r+1 \geq 3k+2t,
    k+t \geq 2,k \geq 1, t \geq 0$.
\par

\section{Acknowledgements} The authors wish to thank R.J.
Gould,  Jiongsheng Li, Rong Luo, J. Schmitt,  Amitabha Tripathi,
Jianhua Yin and Mengxiao Yin for sending their papers to us.

  \par

\end{document}